\documentclass[11pt,a4paper]{article}

\DeclareMathAlphabet{\mathcal}{OMS}{cmsy}{m}{n}
\usepackage[top=25mm,bottom=35mm,left=20mm,right=20mm]{geometry} 

\usepackage{url}
\usepackage{graphicx}
\usepackage{mathptmx}     
\usepackage{mathtools}
\mathtoolsset{showonlyrefs=true}
\DeclarePairedDelimiter\paren{\lparen}{\rparen}
\usepackage{color}
\usepackage{doi}
\usepackage{amsmath}
\usepackage{amsfonts}
\usepackage{amssymb}
\usepackage{amsthm}
\usepackage{bm}

\usepackage[linesnumbered,ruled]{algorithm2e}

\newcommand{\bbR}{\mathbb R}
\newcommand{\tu}{\tilde{u}}
\newcommand{\tv}{\tilde{v}}
\newcommand{\tw}{\tilde{w}}

\newcommand{\tV}{\tilde{V}}
\newcommand{\tW}{\tilde{W}}
\newcommand{\tpsi}{\tilde{\psi}}
\usepackage{bm}
\newcommand{\bb}{\bm b}
\newcommand{\bx}{\bm x}

\newcommand{\bU}{\bm U}
\newcommand{\bV}{\bm V}
\newcommand{\bW}{\bm W}

\newcommand{\rmd}{\mathrm d}
\newcommand{\deltaone}{\delta_x^{\langle 1\rangle}}
\newcommand{\deltatwo}{\delta_x^{\langle 2 \rangle}}
\newcommand{\Dt}{\Delta t}

\theoremstyle{definition}
\newtheorem{scheme}{Scheme}
\newtheorem{proposition}{Proposition}

\usepackage{tikz}
\usepackage{pgfplotstable}
\usetikzlibrary{arrows,positioning,plotmarks,external,patterns,angles,
decorations.pathmorphing,backgrounds,fit,shapes,graphs,calc}

\begin{document}
\title{Modified Strang splitting for semilinear parabolic problems}
\author{Kosuke Nakano\thanks{Department of Applied Physics, Graduate School of Engineering, Nagoya University, Furo-cho, Chikusa-ku, 464-8603 Nagoya,
Japan, k-nakano@na.nuap.nagoya-u.ac.jp}, 
Tomoya Kemmochi\thanks{Department of Applied Physics, Graduate School of Engineering, Nagoya University, Furo-cho, Chikusa-ku, 464-8603 Nagoya,
Japan, kemmochi@na.nuap.nagoya-u.ac.jp},
Yuto Miyatake\thanks{Cybermedia Center, Osaka University, 
			1-32 Machikaneyama, Toyonaka, Osaka 560-0043, Japan,
miyatake@cas.cmc.osaka-u.ac.jp},
Tomohiro Sogabe\thanks{Department of Applied Physics, Graduate School of Engineering, Nagoya University, Furo-cho, Chikusa-ku, 464-8603 Nagoya,
Japan, sogabe@na.nuap.nagoya-u.ac.jp},
Shao-Liang Zhang\thanks{Department of Applied Physics, Graduate School of Engineering, Nagoya University, Furo-cho, Chikusa-ku, 464-8603 Nagoya,
Japan, zhang@na.nuap.nagoya-u.ac.jp}
}
\date{}

\maketitle

\begin{abstract}
We consider applying the Strang splitting to semilinear parabolic problems. The key ingredients of the Strang splitting are the decomposition of the equation into several parts and the computation of approximate solutions by combining the time evolution of each split equation. However, when the  Dirichlet boundary condition is imposed, order reduction could occur due to the incompatibility of the split equations with the boundary condition. In this paper, to overcome the order reduction, a modified Strang splitting procedure is presented for the one-dimensional semilinear parabolic equation with first-order spatial derivatives, like the Burgers equation.
\end{abstract}

\section{Introduction}
\label{sec1}

Splitting methods are numerical methods for solving evolution equations.
The idea, which much differs from other numerical methods such as Runge--Kutta methods, is to decompose the equation into several parts and construct an approximate solution by combining the time integration of each split equation.
There are several splitting methods such as the Lie splitting and Strang splitting.
They differ in the way of combining the solution of each split equation.
In most cases, the Lie splitting is first-order accurate in time, and the Strang is second-order~\cite{mq02}.

In this paper, we are concerned with the numerical integration of the one-dimensional semilinear parabolic equation
\begin{align}
u_t &= u_{xx} + f(u,u_x),  \quad x \in (0,L), \ t > 0,  \label{eq:slp_eq} \\
u (t,0) &=b_1(t), \quad u(t,L) = b_2(t) , \label{slp_bc1}\\
u(0,\cdot) &= u_0,
\end{align}
where $b_1, b_2 :[0,T] \to \bbR$ express the Dirichlet boundary condition.
In this paper, we consider the one-dimensional problem taking in mind that a new approach which will be presented can be generalized to multi-dimensional problems. Further, the new approach remains applicable even if the function $f$ depends on $t$.
We note that discretizing the equation in space leads to a stiff ordinary differential equation system due to the second order derivative, which makes it challenging to apply an explicit Runge--Kutta method for the time discretization.
On the other hand, most implicit methods require solving nonlinear equations if $f$ is nonlinear.
In contrast, it is expected that splitting methods construct an approximate solution more efficiently.
The most straightforward decomposition of the semilinear parabolic equation \eqref{eq:slp_eq} seems the following two equations
\begin{align} \label{naive_splitting}
  v_t = v_{xx} \quad \text{and}\quad w_t = f(w,w_x).
\end{align}
Suppose that these equations are well-posed at least locally.
Let  $\psi_{1,\Dt}$ and $\psi_{2,\Dt}$ denote the time evolution maps for the first and second equations, respectively, i.e. $ \psi_{1,\Dt}: v(t) \mapsto v(t+\Dt) $ and $ \psi_{2,\Dt} : w(t) \mapsto w(t+\Dt)$.
The Strang splitting method~\cite{st68} is formulated as
\begin{align}
	u_{n+1} = \psi_{1,\frac{\Dt}{2}}\circ\psi_{2,\Dt}\circ\psi_{1,\frac{\Dt}{2}}(u_n)
\end{align}
or 
\begin{align}
 	u_{n+1} = \psi_{2,\frac{\Dt}{2}}\circ\psi_{1,\Dt}\circ\psi_{2,\frac{\Dt}{2}}(u_n),   
\end{align}
where $u_n $ denotes the numerical solution at $t=n\Delta t$.
In practice, after appropriate spatial discretization, the $v$-equation is integrated exactly by using the matrix exponential, and the $w$-equation, which is non-stiff in general, by a second-order explicit method such as the Heun method.
As the $v$-equation, which is stiff, is integrated exactly, 
a relatively large time step size can be used. Furthermore, since the $w$-equation is integrated by an explicit method, there is no need to solve nonlinear equations.

However, 
the order of convergence is sometimes deteriorated.
That is, the second-order convergence is not achieved, and this phenomenon is called the order reduction.
Let us consider the case the function $f$ only depends on $u$.
The $v$-equation in \eqref{eq:slp_eq} is nothing but the heat equation, and thus there is no difficulty in defining its time evolution map.
However, the $w$-equation ($w_t=f(w)$) is incompatible with the boundary condition except for the case $\dot{b}_1 = f(b_1) $ and $\dot{b}_2 = f(b_2)$, where the dot stands for the time differentiation.
This incompatibility could cause the order reduction. 
We note that the order reduction also occurs for other type of equations such as the advection-reaction equation~\cite{hv03} (see also~\cite{hv95}).
Einkemmer and Ostermann~\cite{eo15} proposed a new splitting method, which they called the modified Strang splitting, that overcomes the order reduction in multi-dimensional settings (see also~\cite{ac19}).
The key idea is to devise the decomposition after a certain change of variables. 
Convergence analysis was also given in~\cite{eo15}.

In this paper, we consider the case in which the function $f$ also depends on spatial derivatives of $u$ and the order reduction occurs.
We aim to study if the idea proposed by Einkemmer and Ostermann~\cite{eo15} can be extended to the general cases to overcome the order reduction.
In this paper, focusing on the one-dimensional case, in particular Burgers' type, i.e. $f(u,u_x) = uu_x$,
we shall present a modified Strang splitting procedure, and observe that it exhibits the second-order convergence by several numerical experiments.
Discussions for more general cases and theoretical analyses are left to our future work and will be reported elsewhere.
The key idea is
to rewrite the equation as an equation with the homogeneous boundary condition (if the original condition is non-homogeneous) and to devise the decomposition so that the split equations are compatible with the boundary condition.

We use the following notation.
We discretize the space interval $[0,L]$ by an  equidistant grid with the space mesh size $\Delta x = L / (K+1)$.
The approximation to $u(t,k\Delta x)$ is denoted by $U_k(t)$.
We often write the numerical solutions as a vector $\bU(t) = [ U_1(t),\dots, U_K(t)]^\top$.
The first and second order central difference operators are denoted by $\deltaone$ and $\deltatwo$, respectively.
The index for the temporal discretization is denoted by $n$.
For example, $u_n$, $\bU_n$, $U_{n,k}$ are approximations to
$u(n\Dt,\cdot)$, $[u(n\Dt,\Delta x),\dots,u(n\Dt,K\Delta x)]^\top$, $u(n\Delta t, k \Delta x)$, respectively.

\section{Modified Strang splitting}

\subsection{Incompatibility of the simple splitting}
\label{subsec:inc}

As explained in Section~\ref{sec1}, 
the $w$-equation in \eqref{naive_splitting} may be incompatible with the boundary condition.
Below, we explain the incompatibility in more detail.

Let us focus on the case $f(u,u_x) = uu_x$.
The standard decomposition \eqref{naive_splitting} leads to the following two equations
\begin{align}
v_t = v_{xx}\quad \text{and} \quad w_t = ww_x,
\end{align} 
where the boundary condition \eqref{slp_bc1} is imposed to both equations.
We note that the $w$-equation is the advection equation, and it can be uniquely solved by imposing only the left boundary condition $u(t,0)=b_1(t)$.
Therefore, if we impose conditions on both sides of the boundary, the solution to the $w$-equation could be overdetermined unless the boundary condition is homogeneous ($b_1 = b_2 = 0$).
In general, this problem occurs even if the boundary condition is homogeneous, e.g.~$f(u,u_x) = (u+1)u_x$.

We note that 
the situation could change after discretizing the space variable.
Applying the central difference to the $w$-equation gives the semi-discrete scheme $\dot{\bW} = F(\bW)$, where
\begin{align}
F_k (\bW) = W_k \frac{W_{k+1}-W_{k-1}}{2\Delta x}
\end{align}
with $W_0(t) = b_1(t)$ and $W_{K+1}(t) = b_2(t)$.
Because the size of the system and the number of dependent variables coincide,
we can proceed the Strang splitting procedure.
This idea will be tested and compared with a modified Strang splitting which will be given in the next subsection.

\subsection{Modified Strang splitting}

We now present a modified Strang splitting procedure.

First, we rewrite the equation \eqref{eq:slp_eq} 
to avoid the non-homogeneous boundary condition.
For this aim, we introduce a new variable $z$, which is defined by
\begin{align*}
z(t,x) = \frac{b_2(t) - b_1(t)}{L} x + b_1 (t).
\end{align*}
This is the solution to the Laplace equation $z_{xx} (t,\cdot) = 0$ with the boundary condition $z(t,0)=b_1(t)$ and $z(t,L) = b_2(t)$.
Let 
\begin{align}
\tilde{u} := u -z.
\end{align}
Using the variable $\tilde{u}$, we can
rewrite the semilinear parabolic equation \eqref{eq:slp_eq} as the following equation with the homogeneous boundary condition
\begin{align}
\tu_t &= \tu_{xx} + f(\tu + z,\tu_x + z_x) - z_t, \label{eq:vt1} \\
\tu(t,0) &= \tu(t,L) = 0,\\
\tu(0,\cdot) &= u_0 - z(0).
\end{align}

Next, we decompose the equation \eqref{eq:vt1} into two equations.
The simplest decomposition
\begin{align}
    \tv_t = \tv_{xx}  - z_t, \quad
\tw_t = f(\tw + z,\tw_x + z_x)  
\end{align}
does not overcome the incompatibility discussed in Section~\ref{subsec:inc}.
To avoid the incompatibility, we consider the following decomposition
\begin{align}
\tv_t &= \tv_{xx} + f( z,\tv_x + z_x) - z_t, \label{new_sp1}\\
\tw_t &= f(\tw + z,\tw_x + z_x) - f( z,\tw_x + z_x). \label{new_sp2}
\end{align}
The compatibility for the $w$-equation is now achieved, which is intuitively understood in such a way that the both sides of \eqref{new_sp2} are zero at the boundaries.
Detailed discussions are given in the next subsection.
Using the new decomposition, we propose a modified Strang splitting scheme as follows.

\begin{scheme}[Modified Strang splitting]
\label{proposed}
Let  $\tpsi_{1,\Dt}$ and $\tpsi_{2,\Dt}$ denote the time evolution maps for \eqref{new_sp1} and \eqref{new_sp2}, respectively.
The new Strang splitting method is formulated as
\begin{align} \label{mod_Ss1}
	\tu_{n+1} = \tpsi_{1,\frac{\Dt}{2}}\circ\tpsi_{2,\Dt}\circ\tpsi_{1,\frac{\Dt}{2}}(\tu_n)
\end{align}
or
\begin{align}
	\tu_{n+1} = \tpsi_{2,\frac{\Dt}{2}}\circ\tpsi_{1,\Dt}\circ\tpsi_{2,\frac{\Dt}{2}}(\tu_n).
\end{align}
\end{scheme}
In the numerical experiments shown in Section~\ref{sec3}, we will employ \eqref{mod_Ss1}.

\subsection{Compatibility of the $w$-equation}

We explain that the $w$-equation \eqref{new_sp2} is compatible with the boundary condition as long as the initial condition satisfies 
\begin{align}
\tu (0,0) = \tu(0,L) = 0.
\end{align}
For clarity, we focus on the case
$f(u,u_x) = uu_x$.
In this case the $w$-equation is written as 
\begin{align}
\tw_t = \tw (\tw_x + z_x).
\end{align}
To show the compatibility of this equation to the boundary condition, we consider the following problem defined on the whole real line:
\begin{align}
u_t &= uu_x + a(t,x) u, \quad x\in\bbR, \ t>0, \label{eq:onR}\\
u(t,x) &= u_0(x), \quad x \in \bbR.
\end{align}

\begin{proposition}
Assume that there exists $T>0$ such that the solution to \eqref{eq:onR} is sufficiently smooth for $x \in \bbR$ and $t \in [0,T]$.
Then, it follows that $u(t,x_0) =0$ for $t \in [0,T]$ if $u_0(x_0)=0$.
\end{proposition}

Due to the proposition, considering $u=\tw$, 
$a(t,x)=z_x(t,x)$ and $x_0 = 0, L$ indicates the intended compatibility.

\begin{proof}
The proof is based on the method of characteristics.

The characteristic curve $(T(s),X(s))$ of the equation \eqref{eq:onR} is given by 
\begin{align}
\frac{\rmd X(s)}{\rmd s} &=  - u(T(s),X(s)), \\
\frac{\rmd T(s)}{\rmd s} &= 1.
\end{align}
Let us impose the initial conditions $X(0) = x_0$ and $T(0) = 0$.
It is readily obtained that $T(s) = s$.
For the function $h(s) := u(s,X(s))$, it follows that
$h(0) = u(0,X(0)) = u_0(x_0)$ and
\begin{align*}
\frac{\rmd h}{\rmd s} &= -u_x (T(s),X(s)) u(T(s),X(s)) + u_t(T(s),X(s))  \\ &= a(T(s),X(s)) u(T(s),X(s))
= \hat{a}(s)h(s),
\end{align*}
where $\hat{a}(s) = a(T(s),X(s)) $.
Therefore,
\begin{align}
h(s) = u_0(x_0)  \exp \paren*{
\int_0^s \hat{a}(r)\,\rmd r},
\end{align}
which means that
\begin{align}
\frac{\rmd X}{\rmd s} 
= -h(s)
=
-u_0(x_0)  \exp \paren*{
\int_0^s \hat{a}(r)\,\rmd r}.
\end{align}
We thus obtain
\begin{align}
X(s) = x_0 - u_0(x_0)
\int_0^s \exp \paren*{
\int_0^r \hat{a}(q)\,\rmd q
}\,\rmd r.
\end{align}
Hence, if $u_0(x_0) = 0$,
the characteristic curve satisfies $X(s) = x_0$, which indicates that $0=h(s) = u(s,X(s)) = u(s,x_0) $, as long as the solution is sufficiently smooth.
\end{proof}

\section{Implementation}

We discretize the split equations \eqref{new_sp1} and \eqref{new_sp2} by the finite difference method.
We note that the following procedure still makes sense for other discretization such as the finite element method.
We employ the first and second order central differences for $u_x$ and $u_{xx}$.
When $f(u,u_x) = uu_x$, the split equations are discretized as
\begin{align}
\dot{\tV}_k &= \deltatwo \tV_k
+Z_k \deltaone (\tV_k + Z_k) - \dot{Z}_k, \label{semi1} \\
\dot{\tW}_k &= \tW_k \deltaone (\tW_k + Z_k), \label{semi2}
\end{align}
where $Z_k(t) = z(t,k\Delta x)$.

The first equation \eqref{semi1} is written in a compact form as
\begin{align}
\dot{\tilde{\bV}} = A(t) \tilde{\bV} + \bb (t),
\end{align}
where $A(t)$ expresses the linear part, i.e. $A(t)\tilde{\bV}$ corresponds to $(\deltatwo + Z_k \deltaone) \tV_k$ and $\bb(t)$ the remaining part.
For solving this equation, a second-order exponential integrator can be employed~\cite{ho10}.
We note that if the boundary condition is time-independent, $A(t)$ and $\bb(t)$ are constant and the solution
\begin{align}
\tilde{\bV}(t) = \exp (tA)\tilde{\bV}(0) + A^{-1} (\exp (tA)-I) \bb
\end{align}
can be computed exactly by using matrix exponential.
The second equation \eqref{semi2} is non-stiff in general, and thus a second-order explicit method such as the Heun method can be employed. 

We note that the proposed approach requires computing the multiplication of a matrix exponential and vector, i.e. $\phi (tA) \bx = \exp (tA) \bx$.
For readers' convenience, we summarize an efficient way of computing it (see~\cite{ho10,hi08} for more details).
Note that the direct computation of all elements of $\phi (tA)$ is usually prohibitive because it requires all eigenvalues and eigenvectors of the matrix $A$, which may vary as time passes.
All we want to do is to calculate $\exp (tA)\bx$ rather than $\exp(tA)$ itself.
An efficient way of doing this  is to employ a Krylov subspace method.
For example, the Arnoldi procedure generates a matrix $V_m \in \bbR^{n\times m}$ whose column vectors are orthonormal and a Hessenberg matrix $H_m \in \bbR^{m\times m}$ that satisfy
\begin{align}
  V_m^\top (A) V_m = H_m.
\end{align}
Using these matrices, the matrix-vector product can be approximated by
\begin{align}
	\phi(tA) \bx \approx V_m \phi (tH_m) \bx.
\end{align}
In practice, the right-hand side often converges to the left with small $m$.

\section{Numerical experiments}
\label{sec3}

Let us consider the one-dimensional Burgers equation
\begin{align}
  u_t = u_{xx} + uu_x, \label{burgers}
\end{align}
i.e.~$f(u,u_x) = uu_x$.
We check the convergence order numerically for the following two cases.
The reference solution for each case is generated by discretizing the equation in space and integrating the semi-discrete scheme by the explicit fourth-order Runge--Kutta method with the time step size $\Delta t = 10^{-8}$. 

All the computations were performed in a computation environment: 1.6 GHz Intel Core i5, 16GB memory, OS X 10.14. We use Julia version 1.2.0.

\subsection{Case~1}
First, we consider the case the both sides of the boundary take the same value. That is,
we consider the Burgers equation \eqref{burgers}
under the boundary condition
\begin{align}
u(t,0) = u(t,1) = 1,
\end{align}
together with the initial condition
\begin{align}
u(0,x) = 2\sin (\pi x) + 1 , \quad x \in (0,1).
\end{align}

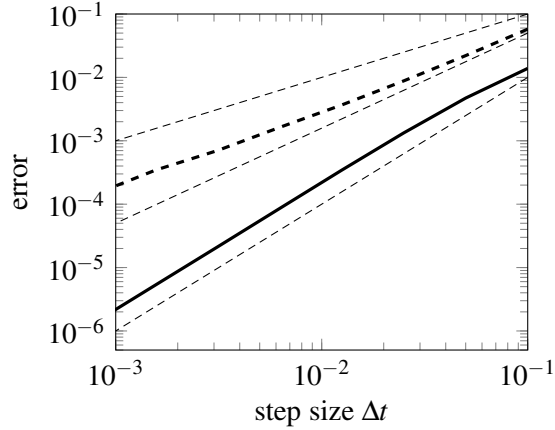
\begin{figure}

\centering
\begin{tikzpicture}
\tikzstyle{every node}=[]
\begin{axis}[width=7cm,
xmax=0.1,xmin=0.001,
ymax=0.1, ymin = 0.0000005,
ymode = log,
xmode = log,
xlabel={step size $\Dt$},ylabel={error},
ylabel near ticks,
legend entries={based on~\eqref{naive_splitting}, Scheme~\ref{proposed}},
legend style={legend cell align=left,draw=none,fill=white,fill opacity=0,text opacity=1,legend pos=south east,},
	]
\addplot[very thick, dashed, 
filter discard warning=false, unbounded coords=discard
] table {
0.1		0.05768862054647905
0.05		0.022071507267696466
0.025		0.008802846715588952
0.0125		0.0036706839196316565
0.00625		0.0016279978152566876
0.003125		0.0007028414530088067
0.0015625		0.00034501156720190274
0.00078125		0.00014117167076821424
};
\addplot[very thick, 
filter discard warning=false, unbounded coords=discard
] table {
0.1		0.013867212348599
0.05		0.00471844909566399
0.025		0.0013298028346901969
0.0125		0.0003436465627320029
0.00625		8.55649602944375e-5
0.003125		2.124875639153423e-5
0.0015625		5.323813066615557e-6
0.00078125		1.331847065744185e-6
};
\addplot[densely dashed,
filter discard warning=false, unbounded coords=discard
] table {
0.1		0.1
0.001		0.001
};
\addplot[densely dashed,
filter discard warning=false, unbounded coords=discard
] table {
0.1		0.05
0.001		0.00005
};
\addplot[densely dashed,
filter discard warning=false, unbounded coords=discard
] table {
0.1		0.01
0.001		0.000001
};
\end{axis}
\end{tikzpicture} 

    \caption{Errors between the numerical solutions and reference solution at $t=0.1$ for Case~1 problem. The thin dashed lines indicate the convergence rates of $1$, $1.5$ and $2$.}
    \label{fig1}
\end{figure}

We compare the numerical behaviour of the modified Strang splitting scheme (Scheme~\ref{proposed}) with that of the Strang splitting based on \eqref{naive_splitting}.
Figure~\ref{fig1} shows errors in the discrete $L^\infty$ norm,
from which
the desired second-order convergence is observed for Scheme~\ref{proposed}.
We note that the error of the Strang splitting based on \eqref{naive_splitting} is much bigger than that of the proposed method, and in particular 
the convergence rate seems deteriorated.

\subsection{Case~2}
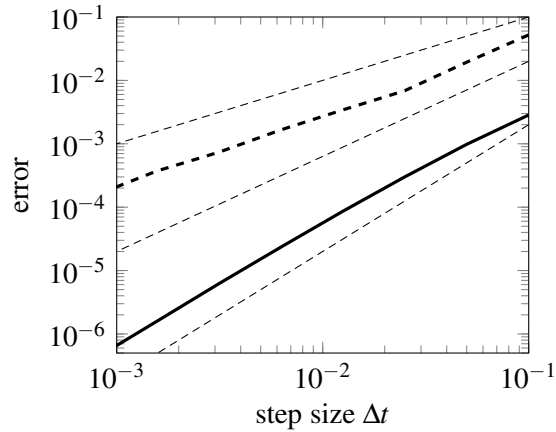
\begin{figure}

\centering
\begin{tikzpicture}
\tikzstyle{every node}=[]
\begin{axis}[width=7cm,
xmax=0.1,xmin=0.001,
ymax=0.1, ymin = 0.0000005,
ymode = log,
xmode = log,
xlabel={step size $\Dt$},ylabel={error},
ylabel near ticks,
legend entries={based on~\eqref{naive_splitting}, Scheme~\ref{proposed}},
legend style={legend cell align=left,draw=none,fill=white,fill opacity=0,text opacity=1,legend pos=south east,},
	]
\addplot[very thick, dashed, 
filter discard warning=false, unbounded coords=discard
] table {
0.1		0.05221804903809213
0.05		0.019381675746031668
0.025		0.00686984511922395
0.0125		0.0034188902071172755
0.00625		0.0016527776578687536
0.003125		0.0007340399131898767
0.0015625		0.00037108189824208715
0.00078125		0.0001527017394258312
};
\addplot[very thick, 
filter discard warning=false, unbounded coords=discard
] table {
0.1		0.00284884121207285
0.05		0.0009806587238141429
0.025		0.000300005775270229
0.0125		8.588044663460082e-5
0.00625		2.3321190898872857e-5
0.003125		6.174526111957235e-6
0.0015625		1.5770268846360125e-6
0.00078125		4.0513243071416127e-7
};
\addplot[densely dashed,
filter discard warning=false, unbounded coords=discard
] table {
0.1		0.1
0.001		0.001
};
\addplot[densely dashed,
filter discard warning=false, unbounded coords=discard
] table {
0.1		0.02
0.001		0.00002
};
\addplot[densely dashed,
filter discard warning=false, unbounded coords=discard
] table {
0.1		0.002
0.001		0.0000002
};
\end{axis}
\end{tikzpicture} 

    \caption{Errors between the numerical solutions and reference solution at $t=0.1$ for Case~2 problem ($b_1=1$, $b_2=3$).}
    \label{fig2}
\end{figure}

We consider the following equation
\begin{align}
& u_t = u_{xx} + uu_x, \quad x \in (0,1), \ t > 0, \\
& u(t,0) =b_1, \quad  u(t,1) = b_2, \\
& u(0,x) = (b_2-b_1)x + b_1 , \quad x \in (0,1).
\end{align}
The results for the case $b_1=1$ and $b_2=3$ are shown in Figure~\ref{fig2}, from which similar behaviour to Case~1 is observed.

\section{Concluding remarks}

In this paper, we have proposed a modified Strang splitting procedure for the semilinear parabolic equations.
Numerical experiments suggested the desired second-order convergence. 

Theoretical analyses with extensions to more general problems such as stochastic differential equations~\cite{kcb17} are left to our future work.
Further, from the computational viewpoint,
a possible drawback of the proposed method is that
when $f(u,u_x)$ nonlinearly depends on $u_x$, the first equation after the decomposition \eqref{new_sp1} becomes nonlinear, and thus its computation is less efficient than linear cases.
It would be interesting to consider if such inefficiency can be avoided.

\section*{Acknowledgments}
This work has been supported in part by JSPS, Japan KAKENHI Grant Numbers 16K17550, 16KT0016, 17H02829, 18H05392, 19K12002, 19K14590 and JST, ACT-I Grant Number JPMJPR18US.

\end{document}